\documentclass[12pt, a4paper, reqno]{amsart}
\usepackage{amscd, amssymb,amsmath,upref}
\allowdisplaybreaks[2]
\begin{document}

\def\C{\mathbb{C}}
\def\Q{\mathbb{Q}}
\def\N{\mathbb{N}}
\def\Z{\mathbb{Z}}
\def\R{\mathbb{R}}
\def\supp{\operatorname{supp}}
\def\Ind{\operatorname{Ind}}
\def\Aut{\operatorname{Aut}}
\def\End{\operatorname{End}}
\newcommand{\Rep}{\operatorname{Rep}}
\def\Ad{\operatorname{Ad}}
\def\sp{\operatorname{span}}
\def\clsp{\overline{\sp}}
\def\T{\mathcal{T}}
\def\dashind{\operatorname{\!-Ind}}
\def\id{\operatorname{id}}
\def\rt{{\operatorname{rt}}}
\def\lt{\operatorname{lt}}

\def\H{\mathcal{H}}
\def\B{\mathcal{B}}

\def\K{\mathcal{K}}
\def\L{\mathcal{L}}

\newcommand{\Chi}{\raisebox{2pt}{\ensuremath{\chi}}}
\newcommand{\under}{\backslash}
\newcommand{\qregular}{regular }

\newcommand{\sknote}[1]{}

\newtheorem{thm}{Theorem}  [section]
\newtheorem{cor}[thm]{Corollary}
\newtheorem{lemma}[thm]{Lemma}
\newtheorem{prop}[thm]{Proposition}
\newtheorem{thm1}{Theorem}
\theoremstyle{definition}
\newtheorem{defn}[thm]{Definition}
\newtheorem{remark}[thm]{Remark}
\newtheorem{example}[thm]{Example}
\newtheorem{remarks}[thm]{Remarks}
\newtheorem{claim}[thm]{Claim}
\newtheorem{problem}[thm]{Problem}
\newtheorem{conj}[thm]{Conjecture}

\numberwithin{equation}{section}

\title[Covariant representations of Hecke algebras]
{Covariant representations of Hecke algebras\\
and imprimitivity for crossed products\\ by homogeneous spaces}

\author[an Huef]{Astrid an Huef}
\address{School of Mathematics\\
The University of New South Wales\\
NSW 2052\\
Australia}
\email{astrid@unsw.edu.au}

\author[Kaliszewski]{S. Kaliszewski}
\address{Department of Mathematics\\Arizona State University\\ AZ
85287-1804\\USA} \email{kaliszewski@asu.edu}

\author[Raeburn]{Iain Raeburn}
\address{School of Mathematical and Physical Sciences\\
University of Newcastle\\
NSW 2308\\ Australia}
\email{Iain.Raeburn@newcastle.edu.au}

\thanks{This research was supported by grants from the Australian Research Council,
the National Science Foundation and the University of New South Wales.}


\date{1 September 2005}

\begin{abstract}
For discrete Hecke pairs $(G,H)$,
we introduce a notion of covariant representation
which reduces in the case where $H$ is normal to the usual definition
of covariance for the action of $G/H$ on $c_0(G/H)$ by
right translation; in many cases where $G$ is a semidirect product,
it can also be expressed in terms of covariance for a semigroup
action.
We use this covariance to characterise the representations of $c_0(G/H)$
which are multiples of the multiplication representation
on $\ell^2(G/H)$,
and more generally, we prove an imprimitivity theorem for
regular representations of certain crossed products by coactions
of homogeneous spaces. We thus obtain new criteria for extending
unitary representations from~$H$ to~$G$.
\end{abstract}

\maketitle

\section*{Introduction}

Let $G$ be a locally compact group, and let $H$ be a closed subgroup
of $G$.
We have recently (\cite{aHRmackey, aHKR-limbo, ahkrw-restrict})
been working on the problem of extending
unitary representations from $H$ to $G$,
using the theory of non-abelian crossed-product
duality.  Our techniques reduce the extension problem
to one of \emph{imprimitivity}; that is,
to deciding whether a certain induced representation
is equivalent to a particular type of regular representation.
At its core, the latter involves characterising the
representations of $C_0(G/H)$ which are equivalent to a
multiple $1\otimes M$ of the representation $M$ by multiplication on $L^2(G/H)$.

If $H$ is normal in $G$, such a 
characterisation can be obtained from the Stone-von Neumann
theorem: a given representation $\nu$ of $C_0(G/H)$
is equivalent to a multiple of $M$ if and only if there
exists a unitary representation $U$ of $G/H$
such that the pair $(\nu,U)$ is covariant
for the action $\rt$ of $G/H$ on $C_0(G/H)$ by right translation.

Motivated by the desire to
extend our techniques to the non-normal case,
in this paper we obtain a similar characterisation which works when
$H$ is a Hecke subgroup of a discrete group $G$.
With the
Hecke algebra $\H(G,H)$ playing the role of the group $G/H$,
we formulate a covariance condition for pairs $(\nu,V)$ of
representations of $c_0(G/H)$ and $\H(G,H)$, and use it to prove
a Stone-von Neumann-type theorem characterising the representations
of $c_0(G/H)$ which are equivalent to a multiple of $M$.
We then use our theorem to solve the imprimitivity problem
mentioned above, and we apply this to obtain
new results on the extension problem for representations.

Since the Hecke algebra does not act on $c_0(G/H)$
in any obvious way, our covariance condition (Definition~\ref{defns})
looks a little unusual.  However, when $H$ is normal in $G$,
representations of $\H(G,H)$ correspond to unitary representations
of $G/H$, and under this correspondence the condition reduces
to the more familiar covariance condition for representations of
$(c_0(G/H),G/H,\rt)$ mentioned above (see Remark~\ref{gp-rem}).
Moreover, in a large class of examples where
$\H(G,H)$ can be realised as a semigroup crossed product,
our covariance condition can be expressed
in terms of the existing notions for group and semigroup actions
(Theorem~\ref{equiv}),
and is also closely related to recent work of
Exel \cite{E} and Larsen \cite{L}
(see Proposition~\ref{exel-prop}).

In somewhat more detail,
our first main result (Theorem~\ref{mult}) states that
the covariant representations of $(c_0(G/H),\H(G,H))$
are precisely those pairs which are equivalent
to a multiple of $(M,\rho)$, where $\rho$ is the natural
representation of $\H(G,H)$ on $\ell^2(G/H)$
analogous to the right-regular representation of $G/H$.
It follows easily that a given representation~$\nu$ of $c_0(G/H)$
is equivalent to $1\otimes M$ if and only if there
exists a representation $V$ of $\H(G,H)$ such that
$(\nu,V)$ is a covariant pair.
Since $G$ is discrete, our proof is more elementary than
that of the Stone-von Neumann theorem, although it is based
on the same observation: covariant pairs
generate sets of operators which behave like matrix units.

In Section~\ref{three}, we use Theorem~\ref{mult} to obtain a new
imprimitivity theorem for $C^*$-crossed products by maximal coactions.
(This is the natural abstract setting for
our application to the extension problem.)
Since $G$ is discrete, a coaction $\delta$ of $G$
on a $C^*$-algebra $B$ is best viewed as a \emph{Fell bundle}
over $G$, which is an analytic version of
a grading of $B$ by~$G$: for each $x\in G$, the set
$B_x = \{ b\in B \mid \delta(b) = b\otimes x \}$ is
a linear subspace of $B$,
we have $B_x B_y \subseteq B_{xy}$ and
$B_x^* = B_{x^{-1}}$ for $x,y\in G$,
and $\cup_{x\in G} B_x$ spans a dense subspace of $B$ (\cite{QuiDC}).
For a subgroup $H$ of $G$, Echterhoff and Quigg (\cite{EQ})
have defined a crossed product $C^*$-algebra $B\times_{\delta|}(G/H)$
which, if $\delta$ is maximal, is universal for suitably covariant pairs of representations
of $B$ and $c_0(G/H)$.
By definition, the \emph{\qregular representations} of this crossed product
are those induced from a representation $\theta$ of $B$ via the covariant
pair $((\theta\otimes\lambda)\circ\delta,1\otimes M)$, where $\lambda$ is the quasi-regular representation of $G$ on $\ell^2(G/H)$.
Our imprimitivity theorem (Theorem~\ref{disc-thm})
characterises these \qregular representations
up to unitary equivalence as those pairs $(\pi,\nu)$
for which there is a representation $V$ of $\H(G,H)$ in the commutant
of $\pi$ such that $(\nu,V)$ is a covariant pair.

When we apply Theorem~\ref{disc-thm} to the extension problem, we
obtain a general result (Theorem~\ref{appltoext})
for representations of $C^*$-dynamical systems involving actions of~$G$.
To see what this says about group representations, recall
that the group $C^*$-algebra $C^*(H)$ is naturally Morita equivalent
to the crossed product $C^*(G)\times_{\delta_G|}(G/H)$, where
$\delta_G$ is
the comultiplication on $C^*(G)$, and thus there is
a bijective correspondence
$U\mapsto (\pi_U, \nu_U)$ between unitary
representations of $H$ and covariant representations
of~$(C^*(G),c_0(G/H))$.  Theorem~1 of
\cite{aHKR-limbo} says that $U$ extends to
a representation of $G$ if and only if the representation of~$C^*(G)\times_{\delta_G|}(G/H)$ corresponding to
$(\pi_U,\nu_U)$ is equivalent to a regular representation.
Thus, by Theorem~\ref{disc-thm}, $U$ extends
if and only if there is a representation $V$ of $\H(G,H)$ in the
commutant of $\pi_U$ such that $(\nu_U,V)$ is a covariant pair.

\subsection*{Conventions}
Let $(G, H)$ be a discrete Hecke pair: this means that
$H$ is a subgroup
of a discrete group $G$ such that every double coset $HxH$
contains just finitely many left cosets.
We use $R$ to denote the right coset counting map,
so that
\[
R(x) = |H\under HxH| = |Hx^{-1}H/H| < \infty
\]
for all $x\in G$.
We view the Hecke algebra
$\H(G, H)$ as the $*$-algebra
of finitely-supported functions on the double-coset space $H\under G/H$
with the operations given
for $HxH\in H\under G/H$ by
\[
fg(HxH) = \sum_{yH\in G/H} f(HyH)g(Hy^{-1}xH)
\qquad\text{and}\qquad
f^*(HxH) = \overline{f(Hx^{-1}H)}.
\]
We write $[HxH]$ for the
characteristic function of the double coset $HxH$, viewed as
an element of $\H(G,H)$, even when $HxH$ happens to be a single
left or right coset.  We use
$\epsilon_{xH}$ to denote the characteristic function of the
coset $xH$, viewed
as an element of
$\ell^2(G/H)$ or  $c_0(G/H)$, and we use $\Chi$ to denote
characteristic functions in other contexts.

All representations of $C^*$-algebras appearing
in this paper are implicitly non-degenerate $*$-homomorphisms.
All representations of Hecke algebras are unital $*$-representations,
and we will often re-state this explicitly for emphasis.

\section{Covariant representations}\label{one}

To understand where our new covariance condition comes from, we
re-examine the group case in more detail.
Let $H$ be a normal subgroup of a locally compact group $G$,
and let $\rho$ be the right-regular representation of $G/H$
on $L^2(G/H)$.
Then the Stone-von Neumann theorem says that the crossed product
$C_0(G/H)\times_\rt(G/H)$ is isomorphic to the
algebra $\K(L^2(G/H))$ of compact
operators via the integrated form $M\times\rho$
of the covariant representation $(M,\rho)$
(see, for example, \cite[Theorem~C.34]{tfb}).
Since every representation of the compacts is equivalent
to a multiple of the identity representation, one
can deduce that every covariant representation
$(\nu,U)$ of $(C_0(G/H),G/H,\rt)$ is equivalent
to a multiple of $(M,\rho)$.

For a discrete group $G$, however, the same conclusion follows
from the more elementary observation that the operators
$\nu(\epsilon_{xH})U(x^{-1}yH)\nu(\epsilon_{yH})$
generate a set of matrix units
as $xH$ and $yH$ run through $G/H$.
It is this approach that we will extend to
Hecke subgroups.

\begin{defn}\label{defns}
Let $(G,H)$ be a discrete Hecke pair, let
$\nu$
be a nondegenerate $*$-representation of $c_0(G/H)$
on a Hilbert space $\H$,
and let $V$ be a unital $*$-representation
of the Hecke algebra $\H(G,H)$ on the same Hilbert space.

We say that $(\nu,V)$ is a \emph{matrix unit pair}
if the collection
\begin{equation*}
\{ \nu(\epsilon_{xH})V([Hx^{-1}yH])\nu(\epsilon_{yH})
\mid xH,yH\in G/H \}
\end{equation*}
is a set of matrix units in $B(\H)$.

We say that $(\nu,V)$ is a \emph{covariant pair}
if
\begin{equation}\label{wc1}
V([HaH])\nu(\epsilon_{xH})V([HbH])
= \sum_{\substack{uH\subseteq Ha^{-1}H\\ vH\subseteq HbH}}
\nu(\epsilon_{xuH})V([Hu^{-1}vH])\nu(\epsilon_{xvH})
\end{equation}
for all $a,x,b\in G$.
\end{defn}

Note that both properties in Definition~\ref{defns}
are preserved by unitary equivalence.

\begin{remark}\label{gp-rem}
When $H$ is normal in $G$, the Hecke algebra is just the
group algebra $\C(G/H)$, and
we can convert between unital $*$-representations
of $\C(G/H)$ and unitary representations of $G/H$
by identifying
group elements in $G/H$ with their characteristic functions.
It follows that a pair $(\nu,V)$ is covariant for $(G,H)$
if and only
if it is covariant in the usual sense for the action $\rt$
of $G/H$ on $c_0(G/H)$ by right translation.  Indeed, in this case
the sums in \eqref{wc1} disappear and we get
\begin{equation}\label{group3}
V(aH)\nu(\epsilon_{xH})V(bH)
= \nu(\epsilon_{xa^{-1}H})V(abH)\nu(\epsilon_{xbH})
\end{equation}
for all $a,x,b\in G$.  Taking $b=a^{-1}$, we recover
the usual covariance condition:
\begin{equation}\label{group2}
V(aH)\nu(\epsilon_{xH})V(aH)^*
= \nu(\epsilon_{xa^{-1}H})
= \nu( \rt_{aH}(\epsilon_{xH})).
\end{equation}
Conversely,
condition~\eqref{group3} can be derived
from~\eqref{group2} by writing
\begin{align*}
V(aH)\nu(\epsilon_{xH})V(bH)
&= V(aH)\nu(\epsilon_{xH})V(aH)^*V(abH)V(b^{-1}H)
\nu(\epsilon_{xH})V(b^{-1}H)^*\\
&= \nu(\rt_{aH}(\epsilon_{xH}))V(abH)
\nu(\rt_{b^{-1}}(\epsilon_{xH}))\\
&=\nu(\epsilon_{xa^{-1}H})V(abH)\nu(\epsilon_{xbH}).
\end{align*}
\end{remark}

It will follow from Proposition~\ref{mu-cov} that every covariant
pair is a matrix unit pair;
we do not know if the converse is true in general
(see also Remark~\ref{mu-rem}).
In the group case, however,
the covariant pairs are exactly the matrix
unit pairs, since then
condition~(iv) of Lemma~\ref{pairs}
is satisfied by every matrix unit pair.

\begin{lemma}\label{pairs}
Let $(G,H)$ be a discrete Hecke pair,
let $\nu$ be a nondegenerate $*$-representation of $c_0(G/H)$,
and let $V$ be a unital $*$-representation of
$\H(G,H)$ on the same Hilbert space $\H$.
Then the following are equivalent\textup:
\begin{enumerate}
\item[(i)]\label{Vnu}
$\displaystyle V([HaH])\nu(\epsilon_{xH})
= \sum_{uH\subseteq Ha^{-1}H}
\nu(\epsilon_{xuH})V([HaH])\nu(\epsilon_{xH})$
for all $a,x\in G$.
\item[(ii)]\label{nuV}
$\displaystyle \nu(\epsilon_{xH})V([HbH])
= \sum_{vH\subseteq HbH}
\nu(\epsilon_{xH})V([HbH])\nu(\epsilon_{xvH})$
for all $x,b\in G$.
\item[(iii)]\label{zero}
$\displaystyle\nu(\epsilon_{xH})V([HbH])\nu(\epsilon_{yH}) = 0$
unless $Hx^{-1}yH=HbH$.
\end{enumerate}
If $(\nu,V)$ is a matrix unit pair, then \textup{(i)}--\textup{(iii)}
are also equivalent to\textup:
\begin{enumerate}
\item[(iv)]\label{norm}
$\bigl\| V([HaH])|_{\nu(\epsilon_{xH})\H} \bigr\|
\leq {R(a)}^{1/2}$
for all $a,x\in G$, where $R$ is the right-coset counting map.
\end{enumerate}
\end{lemma}

\begin{proof}
The equivalence of~(i) and~(ii) is easily seen on taking adjoints.
Suppose condition~(ii) holds.  Then for any $a,x,b\in G$,
\[
\nu(\epsilon_{xH})V([HbH])\nu(\epsilon_{yH})
= \sum_{vH\subseteq HbH}
\nu(\epsilon_{xH})V([HbH])\nu(\epsilon_{xvH})\nu(\epsilon_{yH})
= 0
\]
unless $xtH=yH$ for some $tH\subseteq HbH$, which is precisely
when $x^{-1}yH\subseteq HbH$, which is precisely when
$Hx^{-1}yH=HbH$.  Thus~(ii) implies~(iii).

Next, assume condition~(iii).  Since $\nu$ is nondegenerate,
to establish~(ii) it suffices to show that
\[
\nu(\epsilon_{xH})V([HbH])\nu(\epsilon_{yH})h
= \sum_{vH\subseteq HbH}
\nu(\epsilon_{xH})V([HbH])\nu(\epsilon_{xvH})\nu(\epsilon_{yH})h
\]
for each $y\in G$ and $h\in\H$.  By assumption, the left-hand
side is zero unless $Hx^{-1}yH= HbH$; the right-hand side
is zero unless $xvH=yH$ for some $vH\subseteq HbH$, which is
precisely when $Hx^{-1}yH= HbH$.
When $Hx^{-1}yH\neq HbH$, the sum
on the right collapses to the single term on the left.
Thus~(iii) implies~(ii).

Now suppose $(\nu,V)$ is a matrix unit pair such
that~(i) holds, and fix $a,x\in G$.  Then for each
$uH\subseteq Ha^{-1}H$,  the matrix unit
$\nu(\epsilon_{xuH})V([H(xu)^{-1}xH])\nu(\epsilon_{xH})
=\nu(\epsilon_{xuH})V([HaH])\nu(\epsilon_{xH})$
is a partial isometry with initial projection
$\nu(\epsilon_{xH})$ and final projection
$\nu(\epsilon_{xuH})$.
These final projections are orthogonal for each of
the $R(a)$ different cosets
$uH\subseteq Ha^{-1}H$.
For $h\in\nu(\epsilon_{xH})\H$, we have
$h=\nu(\epsilon_{xH})h$, and hence
\begin{align*}
\| V([HaH])h\|^2
&= \| V([HaH])\nu(\epsilon_{xH})h\|^2\\
&= \Bigl\|\sum_{uH\subseteq Ha^{-1}H}
\nu(\epsilon_{xuH})V([HaH])\nu(\epsilon_{xH})h\Bigr\|^2\\
&= \sum_{uH\subseteq Ha^{-1}H}
\|\nu(\epsilon_{xuH})V([HaH])\nu(\epsilon_{xH})h\|^2\\
&= \sum_{uH\subseteq Ha^{-1}H}
\|h\|^2
= R(a)\|h\|^2.
\end{align*}
Thus~(i) implies~(iv) for matrix unit pairs.

Finally, suppose $(\nu,V)$ is a matrix unit pair such
that~(iv) holds.
Fix $a,x\in G$
and set $P=\sum_{uH\subseteq Ha^{-1}H}\nu(\epsilon_{xuH})$; note that $P$ is a (self-adjoint) projection in $B(\H)$  because the $\nu(\epsilon_{xuH})$'s are mutually orthogonal
projections.
Then for $h\in\H$,
\begin{align*}
R(a)\|\nu(\epsilon_{xH})h\|^2
&\geq \|V([HaH])|_{\nu(\epsilon_{xH})\H}\|^2
\|\nu(\epsilon_{xH})h\|^2\\
&\geq \|V([HaH])\nu(\epsilon_{xH})h\|^2\\
&= \|PV([HaH])\nu(\epsilon_{xH})h\|^2
+ \|(1-P)V([HaH])\nu(\epsilon_{xH})h\|^2\\
&= \Big(\sum_{uH\subseteq Ha^{-1}H}\|\nu(\epsilon_{xuH})
V([HaH])\nu(\epsilon_{xH})h\|^2\Big)
+ \|(1-P)V([HaH])\nu(\epsilon_{xH})h\|^2\\
&= \Big(\sum_{uH\subseteq Ha^{-1}H}\|\nu(\epsilon_{xH})h\|^2\Big)
+ \|(1-P)V([HaH])\nu(\epsilon_{xH})h\|^2\\
\intertext{(since each matrix unit $\nu(\epsilon_{xuH})
V([HaH])\nu(\epsilon_{xH})$ is a partial isometry whose
initial space contains $\nu(\epsilon_{xH})h$)}
&= R(a)\|\nu(\epsilon_{xH})h\|^2 + \|(1-P)V([HaH])\nu(\epsilon_{xH})h\|^2.
\end{align*}
This forces $(1-P)V([HaH])\nu(\epsilon_{xH})h=0$.  Since
$h\in\H$ was arbitrary, this shows that $V([HaH])\nu(\epsilon_{xH})
= PV([HaH])\nu(\epsilon_{xH})$, which is precisely~(i).
This completes the proof.
\end{proof}

\begin{prop}\label{mu-cov}
Let $(G,H)$ be a discrete Hecke pair. The covariant pairs
for $(G,H)$ are precisely the matrix unit pairs which satisfy
the equivalent conditions of Lemma~\ref{pairs}.
\end{prop}

\begin{proof}
For brevity, we set
$\upsilon_{xH,yH}
=\nu(\epsilon_{xH})V([Hx^{-1}yH])\nu(\epsilon_{yH})$,
so $(\nu,V)$ is a matrix unit pair if and only if
$\{ \upsilon_{xH,yH} \mid xH,yH\in G/H\}$
is a set of matrix units.

First suppose $(\nu,V)$ is a covariant pair.
Then for $x, y, z\in G$ we have
\begin{align*}
\upsilon_{xH,yH}\upsilon_{yH,zH}
&=
\nu(\epsilon_{xH})
V([Hx^{-1}yH])\nu(\epsilon_{yH})V([Hy^{-1}zH])
\nu(\epsilon_{zH})\\
&= \nu(\epsilon_{xH})
\Bigl(
\sum_{\substack{uH\subseteq Hy^{-1}xH\\vH\subseteq Hy^{-1}zH}}
\nu(\epsilon_{yuH})V([Hu^{-1}vH])\nu(\epsilon_{yvH})\Bigr)
\nu(\epsilon_{zH}).
\end{align*}
There is exactly one $uH\subseteq Hy^{-1}xH$ such that
$xH=yuH$, namely $uH=y^{-1}xH$; similarly $vH=y^{-1}zH$
is the unique left coset in $Hy^{-1}zH$ such that
$yvH=zH$.  Thus the sums disappear, and
$Hu^{-1}vH = Hx^{-1}yy^{-1}zH = Hx^{-1}zH$, so
the above expression reduces to
$\nu(\epsilon_{xH})V([Hx^{-1}zH])\nu(\epsilon_{zH})
= \upsilon_{xH,zH}$.
Since also
$\upsilon_{xH,yH}^* = \upsilon_{yH,xH}$
and $\upsilon_{xH,yH}\upsilon_{wH,zH}=0$ unless
$yH=wH$ (simply because $\nu$ and $V$ are $*$-homomorphisms),
this shows that $(\nu,V)$ is a matrix unit pair.
Taking $b\in H$ in the covariant pair condition~\eqref{wc1}
shows that $(\nu,V)$ satisfies condition~(i) of Lemma~\ref{pairs}.

Conversely, suppose that $(\nu,V)$ is a matrix unit pair
which satisfies the equivalent conditions of
Lemma~\ref{pairs}.  Note that for any $a,x,b\in G$,
and for any $uH\subseteq Ha^{-1}H$ and $vH\subseteq HbH$,
we have $H(xu)^{-1}xH=HaH$ and $Hx^{-1}(xv)H=HbH$.
Thus, conditions~(i) and~(ii) of Lemma~\ref{pairs}
and the matrix unit assumption give
\begin{align*}
V([HaH])&\nu(\epsilon_{xH})V([HbH])
= V([HaH])\nu(\epsilon_{xH})\, \nu(\epsilon_{xH})V([HbH])\\
&= \sum_{uH\subseteq Ha^{-1}H} \nu(\epsilon_{xuH})
V([HaH])\nu(\epsilon_{xH})\,
   \sum_{vH\subseteq HbH} \nu(\epsilon_{xH})
V([HbH])\nu(\epsilon_{xvH})\\
&= \sum_{\substack{uH\subseteq Ha^{-1}H\\vH\subseteq HbH}}
\nu(\epsilon_{xuH})V([H(xu)^{-1}xH])\nu(\epsilon_{xH})\,
\nu(\epsilon_{xH})V([Hx^{-1}(xv)H])\nu(\epsilon_{xvH})\\
&= \sum_{\substack{uH\subseteq Ha^{-1}H\\vH\subseteq HbH}}
\upsilon_{xuH,xH}\, \upsilon_{xH,xvH}
\quad= \sum_{\substack{uH\subseteq Ha^{-1}H\\vH\subseteq HbH}}
\upsilon_{xuH, xvH}\\
&= \sum_{\substack{uH\subseteq Ha^{-1}H\\vH\subseteq HbH}}
\nu(\epsilon_{xuH})V([Hu^{-1}vH])\nu(\epsilon_{xvH}).
\end{align*}
Hence $(\nu,V)$ is a covariant pair.
\end{proof}

\begin{example}\label{Mrho}
Let $M$ be the representation of $c_0(G/H)$
on $\ell^2(G/H)$ by pointwise multiplication,
and let $\rho$ be the representation of
$\H(G,H)$ on $\ell^2(G/H)$ by right convolution,
so that
\[
M(f)(\epsilon_{yH}) = f(yH)\epsilon_{yH}
\quad\text{and}\quad
\rho([HaH])(\epsilon_{yH})
=\sum_{uH\subseteq Ha^{-1}H} \epsilon_{yuH}
= \Chi_{yHa^{-1}H}
\]
for $f\in c_0(G/H)$ and  $a,y\in G$.
(If $H$ is normal in $G$, then $\rho$ is the representation of the group
algebra $\C(G/H)$ corresponding to the right regular
representation of $G/H$.)
Then $(M,\rho)$ is a covariant pair.
One way to see this is to first compute directly that
for any $a,y\in G$,
\begin{equation}\label{rM}
\rho([HaH])M(\epsilon_{yH})
= \sum_{uH\subseteq Ha^{-1}H}
\epsilon_{yuH}\otimes\overline{\epsilon_{yH}},
\end{equation}
where by definition $\xi\otimes\overline{\eta}(\zeta)
=(\zeta\mid\eta)\xi$
for $\xi,\eta,\zeta\in\ell^2(G/H)$.
Since $M(\epsilon_{xH})=\epsilon_{xH}\otimes\overline{\epsilon_{xH}}$
for any $x\in G$, this gives
\begin{equation}\label{MrM}
M(\epsilon_{xH})\rho([Hx^{-1}yH])M(\epsilon_{yH})
=\sum_{uH\subseteq Hy^{-1}xH}
(\epsilon_{xH}\otimes\overline{\epsilon_{xH}})
(\epsilon_{yuH}\otimes\overline{\epsilon_{yH}})
= \epsilon_{xH}\otimes\overline{\epsilon_{yH}}.
\end{equation}
Thus $(M,\rho)$ is a matrix unit pair which by
\eqref{rM} and \eqref{MrM} satisfies condition~(i) of
Lemma~\ref{pairs}, and hence is a covariant pair
by Proposition~\ref{mu-cov}.
\end{example}

\begin{thm}\label{mult}
Let $(G,H)$ be a discrete Hecke pair,
let $\nu$ be a nondegenerate $*$-representation of $c_0(G/H)$,
and let $V$ be a unital $*$-representation of
$\H(G,H)$ on the same Hilbert space $\H$.  Then\textup:
\begin{enumerate}
\item[(i)]
$(\nu,V)$ is a matrix unit pair if and only if
there exists a Hilbert space $\H_0$
and a representation $\tilde V$ of $\H(G,H)$ on
$\H_0\otimes\ell^2(G/H)$
such that
$(\nu,V)$ is unitarily equivalent to $(1\otimes M,\tilde V)$
and such that
\begin{equation}\label{MR-eq}
\begin{split}
(1\otimes M(\epsilon_{xH}))&\tilde V([Hx^{-1}yH])
   (1\otimes M(\epsilon_{yH}))\\
&= 1\otimes(M(\epsilon_{xH})\rho([Hx^{-1}yH])M(\epsilon_{yH}))
\quad\text{ for all } x,y\in G.
\end{split}
\end{equation}

\item[(ii)]
$(\nu,V)$ is a covariant pair if and only if
there exists a Hilbert space $\H_0$ such that
$(\nu,V)$ is unitarily equivalent to the covariant representation
$(1\otimes M,1\otimes\rho)$ on $\H_0\otimes\ell^2(G/H)$.
\end{enumerate}
In particular, a representation~$\nu$ of $c_0(G/H)$ is equivalent
to a multiple of $M$ if and only if there exists a
representation~$V$ of~$\H(G,H)$
such that $(\nu,V)$ is a covariant pair.
\end{thm}

\begin{proof}
Suppose first that $(\nu,V)$ is a matrix unit pair on $\H$,
and write $\upsilon_{xH,yH}$ for each matrix unit
$\nu(\epsilon_{xH})V([Hx^{-1}yH])\nu(\epsilon_{yH})$.
Now set $\H_0 = \nu(\epsilon_H)\H$.
Then it is straightforward
to verify, using nondegeneracy of $\nu$, that
the rule \[\nu(\epsilon_{zH})\xi \mapsto
\upsilon_{H,zH}\xi\otimes\epsilon_{zH}\quad\quad(\xi\in\H)\]
determines a unitary isomorphism $\Psi$ of
$\H$ onto $\H_0\otimes\ell^2(G/H)$ such that
\[\Ad\Psi(\upsilon_{xH,yH})
= 1\otimes(\epsilon_{xH}\otimes\overline{\epsilon_{yH}})\]
for all $x,y\in G$.
(Equivalently, write
$\H=\bigoplus_{zH\in G/H} \nu(\epsilon_{zH})\H$
and use
$\bigoplus_{zH\in G/H} \upsilon_{H,zH}$
to map $\H$ onto
$\bigoplus_{zH\in G/H} \nu(\epsilon_{H})\H
\cong\H_0\otimes\ell^2(G/H)$.)
In particular, for each $x\in G$, we have
\[
\Ad\Psi(\nu(\epsilon_{xH}))
= \Ad\Psi(\upsilon_{xH,xH})
= 1\otimes(\epsilon_{xH}\otimes\overline{\epsilon_{xH}})
= 1\otimes M(\epsilon_{xH});
\]
thus $\Ad\Psi\circ\nu=1\otimes M$.

Now let $\tilde V = \Ad\Psi\circ V$ be the representation of
$\H(G,H)$ on $\H_0\otimes\ell^2(G/H)$ corresponding to~$V$.
Then for each $x,y\in G$,
\begin{align*}
(1\otimes M(\epsilon_{xH}))\tilde V([Hx^{-1}yH])
(1\otimes M(\epsilon_{yH}))
&= \Ad\Psi(\nu(\epsilon_{xH})V([Hx^{-1}yH])\nu(\epsilon_{yH}))\\
&= \Ad\Psi(\upsilon_{xH,yH})\\
&= 1\otimes(\epsilon_{xH}\otimes\overline{\epsilon_{yH}})\\
&= 1\otimes(M(\epsilon_{xH})\rho([Hx^{-1}yH])M(\epsilon_{yH})),
\end{align*}
using \eqref{MrM} for the last equality.
This proves the forward implication in~(i); the converse
is straightforward because~\eqref{MR-eq} implies that
$(1\otimes M,\tilde V)$ is a matrix unit pair.

Now suppose $(\nu,V)$ is a covariant pair.  Then in particular
$(\nu,V)$ is a matrix unit pair, so we have $\H_0$, $\Psi$,
and $\tilde V$ as in part~(i).  But now
$(1\otimes M,\tilde V)=(\Ad\Psi\circ\nu,\Ad\Psi\circ V)$
is covariant because $(\nu,V)$ is, and it follows that
\[
(1\otimes M(\epsilon_{xH}))
\tilde V([HbH])
(1\otimes M(\epsilon_{yH}))\\
= (1\otimes M(\epsilon_{xH}))
(1\otimes\rho([HbH]))
(1\otimes M(\epsilon_{yH}))
\]
for
\emph{every} $x,b,y\in G$:
this is just \eqref{MR-eq} when
$HbH = Hx^{-1}yH$, and both expressions are zero otherwise
by Proposition~\ref{mu-cov}.
Letting $x$ and $y$ vary now shows that
$\tilde V([HbH])=1\otimes\rho([HbH])$ for each $b$,
so $\tilde V = 1\otimes\rho$.
This proves the forward implication of~(ii), and the
converse is immediate because $(1\otimes M,1\otimes\rho)$
is a covariant pair.

For the last statement of the theorem, it only remains to observe
that if $\Psi$ is a unitary operator intertwining
$\nu$ and $1\otimes M$,
then $V=\Ad\Psi^*\circ(1\otimes\rho)$ is a representation
of $\H(G,H)$ such that $(\nu,V)$ is a covariant pair.
\end{proof}

\begin{remark}\label{mu-rem}
Although it is not necessarily true that every matrix
unit pair is a covariant pair, it follows easily
from Theorem~\ref{mult}
that for
for every matrix unit pair $(\nu,W)$, there exists a
representation $V$ of $\H(G,H)$ such that $(\nu,V)$ is a covariant
pair.
\end{remark}


\section{Hecke algebras which are semigroup crossed products}\label{two}

{}\sknote{I tried to get as much of the technical preliminaries out from
in front of Theorem~\ref{equiv} as I could, both to get the reader to
the statement of the theorem faster and to have those technicalities
appear closer to where they are used.  Most of the stuff is now between
the statement of Theorem~\ref{equiv} and that of Lemma~\ref{key}.
But significantly, I moved most of the discussion of Murphy and Stacey
covariance to after the proof of the theorem.  Also, I moved
the discussion of $\alpha$, and $(\mu,e)$ being Stacey-covariant for $\alpha$,
to just before Corollary~\ref{cor-ext-LL} (formerly a Theorem)
since that's the only place it's really needed.}

Throughout this section,
we consider a cancellative semigroup $S$ satisfying
the Ore condition $Ss\cap St\neq\emptyset$ for all $s,t\in S$,
and we let $Q=S^{-1}S$;
since $S$ is cancellative, the Ore condition implies that
$Q$ is a group containing $S$ as a sub-semigroup (\cite[Theorem~1.24]{CP}).
We suppose that there is an action of $Q$ by automorphisms of a group $N$,
and we let $G=N\rtimes Q$ be the semi-direct product group.
Finally, we let $H$ be a normal subgroup
of $N$ such that $sH=HsH$ for all $s\in S$
and such that there are finitely
many right cosets of $H$ in $HsH$ for each $s\in S$.
Thus $(G,H)$ is a Hecke pair by
\cite[Proposition~1.7]{LL}.%
\footnote{To reconcile our assumption
that $|H\under HsH|<\infty$ with the
assumption that $|s^{-1}Hs\under H|<\infty$ in \cite{LL},
note that for any subgroup $H$ of any group $G$, and for any
$x\in G$, the map $Hxh\mapsto x^{-1}Hxh$ is a bijection
of $H\under HxH$ onto $x^{-1}Hx\under H$.}
The Hecke algebra of Bost and Connes \cite{BC} comes from a
Hecke pair of this sort, with $S=\N^*$, $Q=\Q^*_+$ acting on $N=\Q$
by $q\cdot n=n/q$, and $H=\Z$;
see also \cite[Example~4.3]{brenken} and \cite[Example~2.2]{LL}.

In this situation, the Hecke algebra
$\H(G,H)$ is isomorphic to a certain semigroup crossed product
$\C(N/H)\times_\alpha S$, and thus each representation
$V$ of $\H(G,H)$ corresponds to a suitably covariant pair
$(U,W)$ of representations of $N/H$ and $S$.
The main result of this section (Theorem~\ref{equiv})
shows that the covariance condition~\eqref{wc1} for a pair
$(\nu,V)$ can be expressed
very naturally in terms of
a set of more-familiar covariance conditions
involving only $\nu$, $U$, and $W$.

The precise statement of  Theorem~\ref{equiv} requires
a few preliminaries, which we have endeavored to keep brief;
further discussion of their significance follows the proof of
the theorem.
First, since $s^{-1}Hs\subseteq H$ for $s\in S$,
the formula $(xH)\cdot s=xsH$
defines a right action of $S$ on $G/H$.
We denote by $\rt$ the associated action of $S$ by endomorphisms
of $c_0(G/H)$, so that
$\rt_s(f)(xH)=f(xsH)$ for $s\in S$ and $x\in G$.
We also use $\rt$ to denote
the action of the group $N/H$ by automorphisms of $c_0(G/H)$ induced
from the natural action of $N/H$ on the right of $G/H$.

Next, suppose $\nu$ is a
a nondegenerate representation of $c_0(G/H)$
and $W$ is a representation of $S$ by isometries
on the same Hilbert space.
Motivated by \cite{MurOG},
we say that the pair $(\nu,W)$ is \emph{Murphy-covariant}
for $(c_0(G/H),S,\rt)$ if
\begin{equation*}
W_s\nu(f) = \nu(\rt_s(f))W_s
\quad\text{for all $s\in S$ and $f\in c_0(G/H)$.}
\end{equation*}
(This is a special case of~\eqref{mur-cov} below.)

Lastly, for $s\in S$ and $n\in N$, define elements $\mu_s$ and
$e(nH)$ of $\H(G,H)$ by
\[
\mu_s={R(s)^{-1/2}}[HsH]\quad\text{and}\quad
e(nH)=[HnH],
\]
where $R$ is the right coset counting map.
Then the map $\mu\colon S\to\H(G,H)$ is a representation of $S$
by (algebraic) isometries,
and $e\colon N/H\to\H(G,H)$ is a unitary representation
of $N/H$ (\cite[Theorem~1.9]{LL}).

\begin{thm}\label{equiv}
With $S$, $N$, $G$, and $H$ as described above,
let $\nu$ be a nondegenerate $*$-representation of $c_0(G/H)$,
and let $V$ be a unital $*$-representation of $\H(G,H)$ on
the same space.
Then $(\nu,V)$ is a covariant pair
if and only if
\begin{enumerate}
\item[(i)]
$(\nu,V\circ e)$ is covariant for $(c_0(G/H),N/H,\rt)$, and
\item[(ii)]
$(\nu,V\circ\mu)$ is Murphy-covariant for
$(c_0(G/H),S,\rt)$.
\end{enumerate}
\end{thm}

Observe that~(i) is equivalent to the condition
\begin{equation}\label{rtnh}
V(e(nH))\nu(\epsilon_{xH})V(e(nH))^*
= \nu(\epsilon_{xn^{-1}H})
\quad\text{for all $n\in N$ and $x\in G$}.
\end{equation}
Also, for $x,y\in G$ and $s\in S$ we have
\begin{equation}\label{S-on-coset}
xsH = yH \iff x\in yHs^{-1}=yHs^{-1}H
\iff xH \subseteq yHs^{-1}H.
\end{equation}
(so the endomorphism of $G/H$ determined by $s$ is
precisely $R(s)$-to-one), and
it follows that the endomorphism $\rt_s$
of $c_0(G/H)$
is given on characteristic functions by
\begin{equation*}
\rt_s(\epsilon_{xH})
= \sum_{uH\subseteq Hs^{-1}H} \epsilon_{xuH}.
\end{equation*}
Thus condition~(ii) of Theorem~\ref{equiv}
is equivalent to
\begin{equation}\label{wkly}
V([HsH])\nu(\epsilon_{xH})
= \sum_{uH\subseteq Hs^{-1}H} \nu(\epsilon_{xuH})V([HsH])
\quad\text{for all $s\in S$ and $x\in G$},
\end{equation}
and hence also (on taking adjoints) to
\begin{equation}\label{wkly*}
\nu(\epsilon_{xH})V([Hs^{-1}H])
= \sum_{uH\subseteq Hs^{-1}H} V([Hs^{-1}H])\nu(\epsilon_{xuH})
\quad\text{for all $s\in S$ and $x\in G$}.
\end{equation}

The proof of Theorem~\ref{equiv}
depends on the following lemma, which essentially shows that
certain key operators behave like matrix units.
First recall that for any discrete Hecke pair $(G,H)$,
the product in the Hecke algebra satisfies
\begin{equation}\label{prod}
[HaH][HbH](HxH) = \left| (HaH \cap xHb^{-1}H) / H \right|
\end{equation}
for all $a,x,b\in G$;
if it happens that $HaHbH=HabH$, then
\begin{equation}\label{mult-id}
[HaH][HbH] = \frac{R(a)R(b)}{R(ab)}[HabH]
\end{equation}
by \cite[Lemma~1]{LL-errata}
(see also \cite[Corollary~3.3]{bigdog}).
In our situation, for any $s,t\in S$ and  $n,m\in N$,
taking $a=s^{-1}n$ and $b=mt$
in the above gives the counting formula
\begin{equation}\label{counting}
\left| (Hs^{-1}nH\cap xHt^{-1}m^{-1}H)/H\right|
= \frac{R(s^{-1}n)R(mt)}{R(s^{-1}nmt)}[Hs^{-1}nmtH](HxH)
= \frac{R(t)}{R(s^{-1}t)}
\end{equation}
for any $xH\subseteq Hs^{-1}nmtH$.
(This last identity incorporates the observation
that $R(s^{-1}n)=1$ for all $s\in S$ and $n\in N$,
and more generally that $R(xny)=R(xy)$ for any $x,y\in G$ and $n\in N$.
To see this, note that since $N$ is normal in $G$, $xny=mxy$ for $m=xnx^{-1}\in N$,
and since $H$ is normal in $N$, the rule $Hxyh\mapsto Hmxyh$
gives a well-defined bijection of $H\under HxyH$ onto
$H\under HmxyH = H\under HxnyH$.)

A direct calculation with \eqref{prod}
(as in the proof of \cite[Theorem~1.9]{LL})
shows that for $s\in S$,
\begin{equation}\label{msms}
[HsH][Hs^{-1}H]
= \Chi_{sHs^{-1}}
= \sum_{mH\subseteq sHs^{-1}H}e(mH).
\end{equation}
Also, combining~\eqref{mult-id} and~\eqref{counting} yields
the factorisation
\begin{equation}\label{decomp}
[Hs^{-1}H]e(nH)[HtH]
= \frac{R(t)}{R(s^{-1}t)}[Hs^{-1}ntH]
\end{equation}
for $s,t\in S$ and $n\in N$.
In the important special case where $nH\subseteq tHt^{-1}H$,
so that $HntH=HtH$ and $Ht^{-1}nH=Ht^{-1}H$, we have
\begin{equation}\label{trick}
e(nH)[HtH] = [HtH]
\qquad\text{and}\qquad
[Ht^{-1}H]e(nH)=[Ht^{-1}H].
\end{equation}
This turns out to be the crux of the proof of Lemma~\ref{key}.

\begin{lemma}\label{key}
If $(\nu,V)$ satisfies
the hypotheses and conditions~\textup{(i)}
and~\textup{(ii)} of Theorem~\ref{equiv}, then
\[
V([Hs^{-1}H])\nu(\epsilon_{zH})V(e(kH))V([HrH])
= \nu(\epsilon_{zsH})V([Hs^{-1}krH])\nu(\epsilon_{zkrH})
\]
for each $s,r\in S$, $k\in N$, and $z\in G$.
\end{lemma}

\begin{proof}
By~\eqref{decomp}, then \eqref{wkly} and~\eqref{wkly*}, followed by
\eqref{rtnh},
\begin{align}
\frac{R(r)}{R(s^{-1}r)}
&\nu(\epsilon_{zsH})V([Hs^{-1}krH])\nu(\epsilon_{zkrH})\label{lhs}\\
&=\nu(\epsilon_{zsH})V([Hs^{-1}H])V(e(kH))
V([HrH])\nu(\epsilon_{zkrH})\notag\\
&= \sum_{\substack{uH\subseteq Hs^{-1}H\\vH\subseteq Hr^{-1}H}}
V([Hs^{-1}H])\nu(\epsilon_{zsuH})V(e(kH))
\nu(\epsilon_{zkrvH})V([HrH])\notag\\
&= \sum_{\substack{uH\subseteq Hs^{-1}H\\vH\subseteq Hr^{-1}H}}
V([Hs^{-1}H])\nu(\epsilon_{zsuH})\nu(\epsilon_{zkrvk^{-1}H})
V(e(kH))V([HrH]),\notag\\
\intertext{which, since each term is zero unless
$vH=r^{-1}k^{-1}sukH\subseteq Hr^{-1}H$, reduces to}
&= \sum_{\substack{uH\subseteq Hs^{-1}H\cap\\s^{-1}krHr^{-1}k^{-1}H}}
V([Hs^{-1}H])\nu(\epsilon_{zsuH})V(e(kH))V([HrH]).\label{last}
\end{align}
Each $uH$ in the sum at~\eqref{last} satisfies
$suH\subseteq sHs^{-1}H$
and $k^{-1}sukH\subseteq rHr^{-1}H$,
so repeated use of \eqref{trick}
(together with \eqref{rtnh}) gives
\begin{align*}
V([Hs^{-1}H])&\nu(\epsilon_{zsuH})V(e(kH))V([HrH])\\
&= V([Hs^{-1}H])V(e(suH))\nu(\epsilon_{zsuH})V(e(kH))V([HrH])\\
&= V([Hs^{-1}H])\nu(\epsilon_{zH})V(e(suH))V(e(kH))V([HrH])\\
&= V([Hs^{-1}H])\nu(\epsilon_{zH})V(e(kH))V(e(k^{-1}sukH))V([HrH])\\
&= V([Hs^{-1}H])\nu(\epsilon_{zH})V(e(kH))V([HrH]).
\end{align*}
Thus the terms in the sum at~\eqref{last} are all identical;
by the counting formula~\eqref{counting} there are precisely
$R(r)/R(s^{-1}r)$ of them, so the
sum at~\eqref{last} reduces to
\[
\frac{R(r)}{R(s^{-1}r)}
V([Hs^{-1}H])\nu(\epsilon_{zH})V(e(kH))V([HrH]).
\]
Comparing this with~\eqref{lhs} completes the proof of the lemma.
\end{proof}

\begin{proof}[Proof of Theorem~\ref{equiv}]
First suppose that $(\nu,V)$ is a covariant pair.
To establish~(i) and~(ii), it suffices to check~\eqref{rtnh}
and~\eqref{wkly}, so fix $x\in G$, $n\in N$, and $s\in S$.
For~\eqref{rtnh},
the covariant pair condition
\eqref{wc1} gives
\begin{align*}
V(e(nH))\nu(\epsilon_{xH})V(e(nH))^*
&= V([HnH])\nu(\epsilon_{xH})V([Hn^{-1}H])\\
&= \nu(\epsilon_{xn^{-1}H})V([Hnn^{-1}H])\nu(\epsilon_{xn^{-1}H})\\
&= \nu(\epsilon_{xn^{-1}H}).
\end{align*}
For~\eqref{wkly},
first note that for each
$uH\subseteq Hs^{-1}H$ we have $xusH=xH$ and $Hu^{-1}H=HsH$,
so \eqref{wc1} gives
\[
\nu(\epsilon_{xuH})V([HsH])
= \nu(\epsilon_{xuH})V([HsH])\nu(\epsilon_{xusH})
= \nu(\epsilon_{xuH})V([Hu^{-1}H])\nu(\epsilon_{xH}).
\]
Thus, using \eqref{wc1} again, we have
\begin{align*}
V([HsH])\nu(\epsilon_{xH})
&= \sum_{uH\subseteq Hs^{-1}H}
\nu(\epsilon_{xuH})V([Hu^{-1}H])\nu(\epsilon_{xH})\\
&= \sum_{uH\subseteq Hs^{-1}H}
\nu(\epsilon_{xuH})V([HsH]).
\end{align*}

Conversely,
suppose that $(\nu,V)$ satisfies covariance conditions~(i) and~(ii),
and hence~\eqref{rtnh}, \eqref{wkly}, and~\eqref{wkly*}.
Fix $a,x,b\in G$ and then, using the Ore condition, choose
$s,t,r\in S$ and $n,m\in N$ such that $a=s^{-1}nt$ and $b=t^{-1}mr$.%
\footnote{
First write $a = s_a^{-1}n_at_a$ and $b=t_b^{-1}m_br_b$.
Since $St_a\cap St_b\neq\emptyset$, there exist
$\sigma_a,\sigma_b\in S$ such that $\sigma_at_a=\sigma_bt_b$.
Now set $s=\sigma_as_a$, $n=\sigma_an_a\sigma_a^{-1}$,
$t=\sigma_at_a=\sigma_bt_b$, $m=\sigma_bm_b\sigma_b^{-1}$,
and $r=\sigma_br_b$.}
Then~\eqref{rtnh} and~\eqref{msms} together imply that
for each $pH,qH\subseteq Ht^{-1}H$,
\begin{align*}
\nu(\epsilon_{xpH})V([HtH][Ht^{-1}H])\nu(\epsilon_{xqH})
&= \sum_{mH\subseteq tHt^{-1}H}
\nu(\epsilon_{xpH})V(e(mH))\nu(\epsilon_{xqH})\\
&= \sum_{mH\subseteq tHt^{-1}H}
V(e(mH))\nu(\epsilon_{xpmH})\nu(\epsilon_{xqH})\\
&= V(e(p^{-1}qH))\nu(\epsilon_{xqH}),
\end{align*}
since $mH=p^{-1}qH$ is the unique left coset in $tHt^{-1}H$
such that $xpmH=xqH$.
Using this with~\eqref{wkly} and~\eqref{wkly*} gives
\begin{align*}
V([HtH])\nu(\epsilon_{xH})V([Ht^{-1}H])
&= \sum_{pH,qH\subseteq Ht^{-1}H}
\nu(\epsilon_{xpH})V([HtH][Ht^{-1}H])\nu(\epsilon_{xqH})\\
&= \sum_{pH,qH\subseteq Ht^{-1}H}
V(e(p^{-1}qH))\nu(\epsilon_{xqH}).\\
\end{align*}
Thus, by the factorisation \eqref{decomp} and~\eqref{rtnh} again,
we have
\begin{align}
\frac{R(t)}{R(s^{-1}t)}&\frac{R(r)}{R(t^{-1}r)}
V([HaH])\nu(\epsilon_{xH})V([HbH])\label{lhs2}\\
&= \frac{R(t)}{R(s^{-1}t)}\frac{R(r)}{R(t^{-1}r)}
V([Hs^{-1}ntH])\nu(\epsilon_{xH})V([Ht^{-1}mrH])\notag\\
&= V([Hs^{-1}H])V(e(nH))V([HtH])\nu(\epsilon_{xH})
V([Ht^{-1}H])V(e(mH))V([HrH])\notag\\
&= \sum_{pH,qH\subseteq Ht^{-1}H}
V([Hs^{-1}H])V(e(nH))
V(e(p^{-1}qH))\nu(\epsilon_{xqH})V(e(mH))V([HrH])\notag\\
&= \sum_{pH,qH\subseteq Ht^{-1}H}
V([Hs^{-1}H])\nu(\epsilon_{xpn^{-1}H})V(e(np^{-1}qmH))V([HrH]),\notag
\intertext{which, by Lemma~\ref{key}, becomes}
&=\sum_{pH,qH\subseteq Ht^{-1}H}
\nu(\epsilon_{xpn^{-1}sH})V([Hs^{-1}np^{-1}qmrH])
\nu(\epsilon_{xqmrH}).\label{sum}
\end{align}

The terms in this sum depend only on
$pn^{-1}sH$, which is always contained in $Ht^{-1}n^{-1}sH=Ha^{-1}H$,
and $qmrH$, which is always contained in $Ht^{-1}mrH=HbH$.
Now for each $uH\subseteq Ha^{-1}H$,
by~\eqref{S-on-coset} we have $pn^{-1}sH=uH$
if and only if $pH\subseteq uHs^{-1}nH$,
so by the counting formula~\eqref{counting} the number of left cosets
$pH\subseteq Ht^{-1}H$ such that $pn^{-1}sH=uH$ is
\[
\left|(Ht^{-1}H\cap uHs^{-1}nH)/H\right|
= \left|(Hs^{-1}nH\cap u^{-1}Ht^{-1}H)/H\right|
= \frac{R(t)}{R(s^{-1}t)}.
\]
Similarly, for each $vH\subseteq HbH$,
the number of cosets $qH\subseteq Ht^{-1}H$ such that $qmrH=vH$
is $|(Ht^{-1}H\cap vHr^{-1}m^{-1}H)/H| = R(r)/R(t^{-1}r)$.
Thus we can collect like terms and rewrite~\eqref{sum} as
\[
\frac{R(t)}{R(s^{-1}t)}\frac{R(r)}{R(t^{-1}r)}
\sum_{\substack{uH\subseteq Ha^{-1}H\\ vH\subseteq HbH}}
\nu(\epsilon_{xuH})V([Hu^{-1}vH])\nu(\epsilon_{xvH});
\]
a comparison with~\eqref{lhs2} completes the proof of Theorem~\ref{equiv}.
\end{proof}

\subsection*{Murphy and Stacey covariance}

If $\alpha$ is an action of a semigroup $S$
on a $C^*$-algebra~$A$,
$\pi$ is a nondegenerate representation
$A$ on a Hilbert space $\H$, and $W$ is a representation of
$S$ by isometries of $\H$,
we say that the pair $(\pi,W)$
\emph{Murphy-covariant} for $(A,S,\alpha)$ if
\begin{equation}\label{mur-cov}
W_s\pi(a) = \pi(\alpha_s(a))W_s
\quad\text{ for all $s\in S$ and $a\in A$.}
\end{equation}
We say that $(\pi,W)$
is \emph{Stacey-covariant}
(after \cite{Sta}) if
\begin{equation}\label{sta-cov}
W_s\pi(a)W_s^* = \pi(\alpha_s(a))
\quad\text{ for all $s\in S$ and $a\in A$.}
\end{equation}

Clearly Stacey-covariance implies
Murphy-covariance, but the converse is not true in general.
For instance, it follows from
Example~\ref{Mrho} and Theorem~\ref{equiv} that,
in the context of that theorem,
$(M,\rho\circ\mu)$ is Murphy-covariant
for $(c_0(G/H),S,\rt)$.
Since the endomorphisms $\rt_s$
extend to unital endomorphisms of $c_b(G/H)=M(c_0(G/H))$,
Stacey-covariance would imply that each $\rho(\mu_s)$
is in fact unitary.
However, $\rho(\mu_s)$ cannot be unitary unless $R(s)=1$,
since
by~\eqref{msms}, for each $\epsilon_{xH}\in\ell^2(G/H)$ we have
\begin{align*}
\rho(\mu_s)\rho(\mu_s)^*(\epsilon_{xH})
&= \frac1{R(s)}\rho([HsH][Hs^{-1}H])(\epsilon_{xH})\\
&= \frac1{R(s)}\sum_{mH\subseteq sHs^{-1}H}
\rho(e(mH))(\epsilon_{xH})\\
&= \frac1{R(s)}\sum_{mH\subseteq sHs^{-1}H}
\epsilon_{xm^{-1}H}
= \frac1{R(s)}\Chi_{xsHs^{-1}H}.
\end{align*}

\subsection*{Exel covariance}

In recent work of Exel~(\cite{E}),
a crossed product $A\times_{\alpha, L}\N$ is associated to
each pair $(\alpha,L)$ consisting of an endomorphism $\alpha$ of a $C^*$-algebra $A$ and a ``{transfer function}'' $L$ for $\alpha$.
In \cite{BR}
(and implicitly in \cite{E}), it was shown that Exel's
crossed product is universal for a family of representations which are,
among other things,  Murphy-covariant for the semigroup
action $(A,\N,\alpha)$.
As we will show below, this family includes representations which
arise naturally from covariant pairs.%

In the context of Theorem~\ref{equiv},
fix $s\in S$ and define $L_s\colon c_0(G/H)\to c_0(G/H)$ by
\begin{equation}\label{defL}
L_s(f)(xH)=\frac{1}{R(s)}\sum_{yH\in G/H\colon ysH=xH}f(yH).
\end{equation}
Then $L_s$ is a positive bounded linear map, and it is a
\emph{transfer function} for the endomorphism $\rt_s$ in the sense that $L_s(\rt_s(g)f)=gL_s(f)$ for $f,g\in c_0(G/H)$,
since for $xH\in G/H$ we have
\[
L_s(\rt_s(g)f)(xH)
=\frac{1}{R(s)}\sum_{yH\in G/H\colon ysH=xH}g(ysH)f(yH)
=g(xH)L_s(f)(xH).
\]
In following proof,
we will frequently make use of the easily-verified formula
\begin{equation}\label{Lonepsilon}
L_s(\epsilon_{zH})=\frac{1}{R(s)}\epsilon_{zsH}.
\end{equation}

\begin{prop}\label{exel-prop}
Let $(\nu,V)$ be a covariant pair for $(G,H)$,
with $G$, $H$, and $S$ as in Theorem~\ref{equiv}.
Then for each $s\in S$,
the pair $(\nu, V(\mu_s))$ is a covariant representation of $(c_0(G/H),\rt_s,L_s)$ in the sense of \cite[Definition~3.8]{BR}.
\end{prop}

\begin{proof}
We need to verify conditions~(TC1) and~(TC2)
of~\cite[Definition~3.1]{BR}, and condition~(C3)
of~\cite[Definition~3.8]{BR}.
The first of these follows immediately from
the Murphy-covariance of $(\nu,V\circ\mu)$ in Theorem~\ref{equiv}.
For each $zH\in G/H$,
taking $r=s$ and $k=e$ in Lemma~\ref{key}
and invoking \eqref{Lonepsilon} gives
\begin{align*}
V(\mu_s)^*\nu(\epsilon_{zH})V(\mu_s)
&=\frac{1}{R(s)}V([Hs^{-1}H])\nu(\epsilon_{zH})V([HsH])\\
&=\frac{1}{R(s)}\nu(\epsilon_{zsH})
=\nu(L_s(\epsilon_{zH})),
\end{align*}
which implies~(TC2).
(Thus $(\nu, V(\mu_s))$ is
a \emph{Toeplitz-covariant} representation of $(c_0(G/H),\rt_s,L_s)$).

To verify the remaining covariance condition~(C3), we need to know%
\sknote{Actually, whatever $K_{\rt_s}$ is, it's contained
in $c_0(G/H)$, and we end up proving that \eqref{Exelcov}
holds for all $f\in c_0(G/H)$.  So we don't really need to know
what $K_{\rt_s}$ is. But we do need the calculation
of $\phi_s(\epsilon_{zH})$.}
which
elements $f\in c_0(G/H)$ act as compact operators on the right-Hilbert
$c_0(G/H)$-module $M_{L_s}$ obtained by completing $c_0(G/H)$ in the norm
defined by the $c_0(G/H)$-valued inner product
$\langle f,g\rangle_{L_s} = L_s(f^*g)$ (see the start of \cite[\S3]{BR}).
Write $\phi_s$  for the
homomorphism of $c_0(G/H)$ into $\L(M_{L_s})$ defined by the left action of
$c_0(G/H)$, which is given by ordinary multiplication on $c_0(G/H)\subseteq
M_{L_s}$. We claim that $\phi_s(\epsilon_{zH})$ is the rank-one operator
$R(s)\Theta_{\epsilon_{zH},\epsilon_{zH}}$,
where by definition $\Theta_{g,h}(f)=g\langle h,f\rangle_{L_s}$.
Indeed, for $f\in c_0(G/H)\subseteq M_{L_s}$ and $xH\in G/H$
we have, using (\ref{Lonepsilon}) again,
\begin{align*}
R(s)\Theta_{\epsilon_{zH},\epsilon_{zH}}(f)(xH)
&=R(s)(\epsilon_{zH}\langle \epsilon_{zH},f\rangle)_{L_s}(xH)\\
&=R(s)\epsilon_{zH}(xH)\rt_s(L_s(\epsilon_{zH}^*f))(xH)\\
&=R(s)\epsilon_{zH}(xH)\rt_s(f(zH)L_s(\epsilon_{zH}))(xH)\\
&=R(s)\epsilon_{zH}(xH)f(zH)L_s(\epsilon_{zH})(xsH)\\
&=\epsilon_{zH}(xH)f(zH)\epsilon_{zsH}(xsH)\\
&=\epsilon_{zH}(xH)f(xH),
\end{align*}
which is just $(\phi_s(\epsilon_{zH})f)(xH)$, proving the claim.
It follows that $\phi_s$ takes values in $\K(M_{L_s})$,
so the ideal $K_{\rt_s}$ appearing in \cite[Definition~3.8]{BR}
is all of $c_0(G/H)$.

The covariance condition~(C3) is expressed in terms of the
representation $\psi_s$
of $M_{L_s}$ defined on the dense subspace $c_0(G/H)$
by $\psi_s(f)=\nu(f)V(\mu_s)$.
($\psi_s$ would be denoted by $\psi_{V(\mu_s)}$ in \cite{BR}.)
It requires that
\begin{equation}\label{Exelcov}
(\psi_s,\nu)^{(1)}(\phi_s(f))=\nu(f)\quad \mbox{for $f\in K_{\rt_s}$},
\end{equation}
where $(\psi_s,\nu)^{(1)}$ is the representation of $\K(M_{L_s})$
induced by the Toeplitz representation   $(\psi_s,\nu)$ of $M_{L_s}$
(see \cite[\S2 and Lemma~3.2]{BR}).
By continuity
it suffices to check \eqref{Exelcov} for $f=\epsilon_{zH}$.
Using~\eqref{msms}, we have
\begin{align*}
(\psi_s,\nu)^{(1)}(\phi_s(\epsilon_{zH}))
&=R(s)(\psi_s,\nu)^{(1)}(\Theta_{\epsilon_{zH},\epsilon_{zH}})\\
&=R(s)\psi_s(\epsilon_{zH})\psi_s(\epsilon_{zH})^*\\
&=R(s)\nu(\epsilon_{zH})V(\mu_s)V(\mu_s)^*\nu(\epsilon_{zH})\\
&=\nu(\epsilon_{zH})V([HsH][Hs^{-1}H])\nu(\epsilon_{zH})\\
&=\sum_{mH\subseteq sHs^{-1}H}\nu(\epsilon_{zH})V(e(mH))\nu(\epsilon_{zH}),
\end{align*}
but by covariance property~(i) in Theorem~\ref{equiv},
\[
\nu(\epsilon_{zH})V(e(mH))\nu(\epsilon_{zH})
=V(e(mH))\nu(\epsilon_{zmH})\nu(\epsilon_{zH})
=\begin{cases}\nu(\epsilon_{zH})&\text{if $mH=H$}\\
0&\text{otherwise.}\end{cases}
\]
Thus the above sum collapses to the single term
$\nu(\epsilon_{zH})$,
as required.
\end{proof}

\begin{remark}
One can verify directly that $R(s)R(t)=R(ts)$ for $s,t\in S$.
It then follows from another application of~\eqref{Lonepsilon}
that $L_sL_t=L_{ts}$,
so $(\alpha,L)$ is an action of~$S$ on $c_0(G/H)$
in the sense of~Larsen \cite{L}.
(The extra condition on extendibility of~$L_s$ required
in~\cite{L} holds because the formula \eqref{defL} also defines an
operator $L_s$ on $c_b(G/H)=M(c_0(G/H))$.)
\end{remark}


\section{The imprimitivity theorem}\label{three}

{}\sknote{I've streamlined (a bit) the proof of the converse direction
of the proof of Theorem~\ref{disc-thm} to fully use the assumption
that $(\nu, V)$ is a covariant pair, not just a matrix unit pair.
This is partly because I found it confusing to read the proof the
old way --- it's not explained until the following remark
(about there existing a cov't pair whenever there exists a mu pair)
why we
didn't use the full strength of the hypotheses.  Also, though,
the main point of that remark didn't have to do with $B$, and
it's an easy direct consequence of Theorem~\ref{mult}, so I've moved
it there.  Finally, the old remark was a bit misleading: we get
$\pi\times\nu$ is equiv to a reg rep'n if and only if there exists
$W$ such that $(\nu,W)$ is a matrix unit pair
\emph{and the ranges of $W$ and $\pi$ commute}.  Taking
$V=\Ad\Psi^*\circ(1\otimes\rho)$ now gives a covariant pair
$(\nu,V)$, but there's no reason that the ranges of $\pi$
and $V$ should commute.  }

Following Echterhoff and Quigg in \cite{EQ},
let $\delta\colon B\to B\otimes C^*(G)$ be a
coaction of a discrete group $G$ on a $C^*$-algebra $B$,
and let $H$ be a subgroup of $G$. The
spectral subspaces $B_x=\{b\in B\mid\delta(b)=b\otimes x\}$ form a
Fell bundle $\B$ over $G$,  and the direct product $\B\times
G/H$ is a Fell bundle over the transformation groupoid
$G\times G/H$. The $*$-algebra
$\Gamma_c(\B\times G/H)$ of finitely supported sections
of this bundle,
which we identify with $\sp\{ (b,xH) \mid b\in B, xH\in G/H\}$,
has an (universal) enveloping $C^*$-algebra which is
denoted by $C^*(\B\times G/H)$ in \cite{EQ}.

Recall from \cite[Definition~2.4]{EQ} that
representations $\pi$ of $B$ and  $\nu$ of
$c_0(G/H)$ on the same Hilbert space form a
\emph{covariant representation of
$(B,G/H, \delta|)$} if
\begin{equation}\label{defn-cov-rep}
\pi(b_x)\nu(\epsilon_{yH})=\nu(\epsilon_{xyH})\pi(b_x)\quad\text{for all $x,y\in
G, b_x\in B_x$.}
\end{equation}
(Note that when $H$ is not normal in $G$, the notation $\delta|$ is
purely formal: it does not stand for an actual coaction.)
For every covariant representation $(\pi,\nu)$ of $(B,G/H,
\delta|)$ there is a unique  representation $\pi\times\nu$ of
$C^*(\B\times G/H)$, called the \emph{integrated form} of $(\pi,\nu)$,
such that
$\pi\times\nu(b,xH)=\pi(b)\nu(\epsilon_{xH})$
for all $b\in B$ and $x\in G$.
We shall assume that $\delta$ is maximal in the sense of
\cite{EKQ}; since  $G$ is discrete this is equivalent to
asking that the $C^*$-algebra $B$ be isomorphic to $C^*(\B)$,
the enveloping $C^*$-algebra of $\Gamma_c(\B)$ (\cite[Proposition~4.2]{EKQ}).
In this case, we know from \cite[Proposition~2.7]{EQ}
that every representation of $C^*(\B\times G/H)$
comes from a covariant pair.
This universal property implies in particular
that when $N$ is a normal subgroup of $G$ ---
so that the restriction $\delta|$ is a coaction of the group $G/N$ ---
$C^*(\B\times G/N)$ is isomorphic to the usual crossed
product $B\times_{\delta|}(G/N)$
(\cite[Corollary~2.8]{EQ}).
Thus for general $H$ we feel free to write
$B\times_{\delta|}(G/H)$ for $C^*(\B\times G/H)$,
which we call the \emph{Echterhoff-Quigg crossed product}.

\begin{example}\label{quasi-ex}
Let $\theta$ be a representation of $B$ on a Hilbert space $\H_\theta$,
and denote by
$\lambda$ the quasi-regular representation of $G$ on $\ell^2(G/H)$.
Then the pair $((\theta\otimes\lambda)\circ\delta, 1\otimes M)$ is a
covariant representation of $(B,G/H, \delta|)$ on
$\H_\theta\otimes\ell^2(G/H)$.
For $x,y\in G$, 
$\lambda_x(\epsilon_{yH}) = \epsilon_{xyH}$,
so  for  $b_x\in B_x$ we have
\begin{align*}
(\theta\otimes\lambda)\circ\delta(b_x)(1\otimes M)(\epsilon_{yH})
&=(\theta(b_x)\otimes \lambda_x)(1\otimes M(\epsilon_{yH}))\\
&=\theta(b_x)\otimes (\lambda_xM(\epsilon_{yH})\lambda_x^*\lambda_x)\\
&=(1\otimes M(\epsilon_{xyH}))(\theta(b_x)\otimes \lambda_x)\\
&=(1\otimes M)(\epsilon_{xyH})(\theta\otimes\lambda)\circ\delta(b_x).
\end{align*}
\end{example}

Just as for normal subgroups,
we call the integrated form $((\theta\otimes\lambda)\circ\delta)
\times(1\otimes M)$ of the covariant representation in
Example~\ref{quasi-ex} the
\emph{\qregular representation of $B\times_{\delta|}(G/H)$
induced by $\theta$}.
The following imprimitivity theorem characterises these
representations for Hecke subgroups.

\begin{thm}\label{disc-thm}
Let $\delta$ be a maximal  coaction of a discrete group $G$ on a
$C^*$-algebra $B$, and let $H$ be a Hecke subgroup of $G$.  Suppose
$\pi\times\nu$ is a representation of
the Echterhoff-Quigg crossed product
$B\times_{\delta|}(G/H)$ on a Hilbert
space $\H$. Then
$\pi\times\nu$ is equivalent to a \qregular representation
if and only if
there exists a representation $V$ of $\H(G,H)$ on $\H$
such that $(\nu,V)$ is a covariant pair and such that
the range of $V$ commutes with the range of $\pi$.
\end{thm}

\begin{proof}
First suppose that $\pi\times\nu$ is equivalent to a \qregular
representation, so
there exists a representation $\theta\colon B\to
B(\H_\theta)$ and a unitary isomorphism
$\Psi\colon\H\to \H_\theta\otimes\ell^2(G/H)$ which intertwines
$\pi\times\nu$ and $((\theta\otimes\lambda)\circ\delta)
\times(1\otimes M)$.
Set $V=\Ad\Psi^*\circ(1\otimes\rho)$;
then $(\nu,V)$ is a covariant pair because it is equivalent
to $(1\otimes M,1\otimes\rho)$,
and the ranges
of $\pi$ and $V$ commute because the ranges of
$(\theta\otimes\lambda)\circ\delta$ and $1\otimes\rho$ commute.
This proves the forward implication.

For the converse, suppose there exists $V$ such that
 $(\nu,V)$ is a covariant pair and such that the
ranges of $\pi$ and $V$ commute.  Then by Theorem~\ref{mult},
there exists a Hilbert space $\H_0$
and a unitary isomorphism
$\Psi\colon\H\to\H_0\otimes\ell^2(G/H)$
which intertwines $(\nu,V)$ and
$(1\otimes M,1\otimes\rho)$.
Set $\tilde\pi = \Ad\Psi\circ\pi$, and
for conciseness, also set
$\tilde\nu = \Ad\Psi\circ\nu = 1\otimes M$
and $\tilde V = \Ad\Psi\circ V = 1\otimes\rho$.

We need to produce a representation
$\theta\colon B\to B(\H_0)$ such that
$(\theta\otimes\lambda)\circ\delta = \tilde\pi$.
To do this, first fix $z\in G$ and $b_z\in B_z$.
We will show that $(1\otimes\lambda_z^*)\tilde\pi(b_z)$
commutes with $1\otimes\K(\ell^2(G/H))$;
it will then follow that there exists $\theta(b_z)\in B(\H_0)$
such that
$\theta(b_z)\otimes 1
=(1\otimes\lambda_z^*)\tilde\pi(b_z)$.
It suffices to check that $(1\otimes\lambda_z^*)\tilde\pi(b_z)$
commutes with the operators (recall~\eqref{MrM})
\[
1\otimes(\epsilon_{xH}\otimes\overline{\epsilon_{yH}})
= 1\otimes(M(\epsilon_{xH})\rho([Hx^{-1}yH])M(\epsilon_{yH}))
= \tilde\nu(\epsilon_{xH})\tilde V([Hx^{-1}yH])\tilde\nu(\epsilon_{yH})
\]
for $x,y\in G$.
Since $(\tilde\pi,\tilde\nu)$ is covariant for $(B,G/H,\delta|)$,
and since the ranges of $\tilde\pi$ and $\tilde V$ commute,
we have
\begin{align*}
(1\otimes\lambda_z^*)\tilde\pi(b_z)
(1\otimes(\epsilon_{xH}\otimes\overline{\epsilon_{yH}}))
&=
(1\otimes\lambda_z^*)\tilde\pi(b_z)\tilde\nu(\epsilon_{xH})
\tilde V([Hx^{-1}yH])\tilde\nu(\epsilon_{yH})\\
&=
(1\otimes\lambda_z^*)\tilde\nu(\epsilon_{zxH})\tilde\pi(b_z)
\tilde V([Hx^{-1}yH]) \tilde\nu(\epsilon_{yH})\\
&=
(1\otimes\lambda_z^*)\tilde\nu(\epsilon_{zxH})
\tilde V([Hx^{-1}yH])\tilde\pi(b_z) \tilde\nu(\epsilon_{yH})\\
&=
(1\otimes\lambda_z^*)\tilde\nu(\epsilon_{zxH})\tilde V([Hx^{-1}yH])
\tilde\nu(\epsilon_{zyH})\tilde\pi(b_z).\\
\intertext{By the covariance of
 $(\tilde\nu, 1\otimes\lambda)=(1\otimes M,1\otimes\lambda)$
for the action of $G$ by left translation on $c_0(G/H)$, and
because the
ranges of $\tilde V=1\otimes\rho$ and $1\otimes\lambda$ commute,
this is}
&=
\tilde\nu(\epsilon_{xH})(1\otimes\lambda_z^*)\tilde V([Hx^{-1}yH])
\tilde\nu(\epsilon_{zyH})\tilde\pi(b_z)\\
&=
\tilde\nu(\epsilon_{xH})\tilde V([Hx^{-1}yH])
(1\otimes\lambda_z^*)\tilde\nu(\epsilon_{zyH})\tilde\pi(b_z)\\
&=
\tilde\nu(\epsilon_{xH})\tilde V([Hx^{-1}yH])
\tilde\nu(\epsilon_{yH})
(1\otimes\lambda_z^*)\tilde\pi(b_z)\\
&=
(1\otimes(\epsilon_{xH}\otimes\overline{\epsilon_{yH}}))
(1\otimes\lambda_z^*)\tilde\pi(b_z).
\end{align*}

Now linearly extend the assignment $b_z\mapsto\theta(b_z)$
to get a map $\theta$ on
$\Gamma_c(\B)=\sp\{ b_z \mid z\in G, b_z\in B_z\}$
which, by construction, satisfies
$(\theta\otimes\lambda)\circ\delta(b) = \tilde\pi(b)$
for $b\in \Gamma_c(\B)$.
Then $\theta$ is in fact a $*$-homomorphism, since
for $b_z\in B_z$ (so that $(b_z)^*\in B_{z^{-1}}$)
we have
\begin{align*}
\theta((b_z)^*)\otimes 1
&= (1\otimes\lambda_{z^{-1}}^*)\tilde\pi((b_{z})^*)
= \tilde\pi((b_{z})^*)(1\otimes\lambda_{z^{-1}}^*)\\
&= \tilde\pi(b_{z})^*(1\otimes\lambda_{z})
= ((1\otimes\lambda_{z}^*)\tilde\pi(b_{z}))^*
= \theta(b_z)^*\otimes 1,
\end{align*}
and in addition, for $w\in G$ and $c_w\in B_w$ we have
\begin{align*}
\theta(b_z)\theta(c_w)\otimes 1
&=(\theta(b_z)\otimes 1)(1\otimes\lambda_w^*)\tilde\pi(c_w)
=(1\otimes\lambda_w^*)(\theta(b_z)\otimes 1)\tilde\pi(c_w)\\
&=(1\otimes\lambda_w^*)(1\otimes\lambda_z^*)\tilde\pi(b_z)\tilde\pi(c_w)
=(1\otimes\lambda_{zw}^*)\tilde\pi(b_zc_w)
=\theta(b_zc_w)\otimes 1.
\end{align*}
By the universal property of $C^*(\B)$, $\theta$ therefore
extends to a representation of $B\cong C^*(\B)$
such that $(\theta\otimes\lambda)\circ\delta = \tilde\pi$,
as desired.
\end{proof}

Suppose, as in Section~\ref{two}, that $G$ is a semidirect product
$N\rtimes Q$, where $Q=S^{-1}S$
for an Ore semigroup~$S$, and  that $H$ is a normal subgroup of $N$ such that
$sH=HsH$ and $|H\under HsH|<\infty$ for all $s\in S$.
Then Theorem~\ref{disc-thm} can be re-formulated without reference to the Hecke
algebra.  To do so, we first recall from
Theorem~1.9 and Corollary~1.12 of \cite{LL}
(see also \cite{LL-errata}) that
there is an action $\alpha$ of $S$
by injective corner endomorphisms of $\C(N/H)$
such that the pair $(e,\mu)$ appearing in Theorem~\ref{equiv}
gives rise to an isomorphism $e\times\mu$
of the $*$-algebraic semigroup crossed
product $\C(N/H)\times_\alpha S$ onto $\H(G,H)$.

Now, blurring the distinction between representations of
$N/H$, $\C(N/H)$, and $C^*(N/H)$, it follows from the
universal property of the crossed product
that the map $V\mapsto (V\circ e,V\circ\mu)$
is a bijection between the set of unital
$*$-representations of $\H(G,H)$ and
the  set of Stacey-covariant pairs (see~\eqref{sta-cov})
for $(C^*(N/H),S,\alpha)$.
Moreover, the range of a given representation $V$ of $\H(G,H)$
is precisely the
$*$-algebra generated by $V(e(N/H))$ and $V(\mu(S))$.
The following corollary is now immediate from Theorem~\ref{disc-thm}.

\begin{cor}\label{cor-ext-LL}
Let $S$, $N$, $G$, and $H$ be as above,
let $\delta$ be a maximal coaction of $G$ on a $C^*$-algebra $B$,
and let $\pi\times\nu$ be a
representation of the Echterhoff-Quigg crossed product
$B\times_{\delta|}(G/H)$ on
a Hilbert space $\H$.
Then $\pi\times\nu$ is equivalent to a \qregular representation
if and only if
there is a unitary representation $U\colon N/H\to B(\H)$ and an isometric
representation $W\colon S\to B(\H)$ such that
\begin{enumerate}
\item[(i)]
$(U,W)$ is a Stacey-covariant representation of $(C^*(N/H),S,\alpha)$,
\item[(ii)]
$U(N/H)$ and $W(S)$ commute with the range of $\pi$,
\item[(iii)]
$(\nu,W)$ is a Murphy-covariant representation of
$(c_0(G/H),S,\rt)$, and
\item[(iv)]
$(\nu, U)$ is a covariant representation of $(c_0(G/H), N/H,\rt)$.
\end{enumerate}
\end{cor}

It is interesting to see here a natural juxtaposition of
Murphy-covariance and Stacey-covariance for the same semigroup.

\section{Application to the extension of unitary representations}
\label{four}

To apply our imprimitivity theorem to the extension problem
considered in  \cite{aHKR-limbo}, we want to take $(B,\delta)$ to
be a dual coaction $(A\times_\alpha G,\widehat\alpha)$. We know
from \cite[Proposition~3.4]{EKQ}
that $\widehat\alpha$ is maximal, so
Theorem~\ref{disc-thm} applies, but we also need to know that the
Echterhoff-Quigg
crossed product $(A\times_\alpha G)\times_{\widehat\alpha|}
(G/H)$ is isomorphic to the imprimitivity algebra $(A\otimes
c_0(G/H))\times_{\alpha\otimes\lt}G$ which appears in
\cite[Theorem~1]{aHKR-limbo}.
(Throughout this section,
$\lt$ denotes the action of $G$ on $c_0(G/H)$ by left translation.)
Related, but different, versions of the following proposition
can be found in \cite[Lemma~2.4]{EKR}
and \cite[Theorem~A.64]{EKQR}
(for normal subgroups),
and \cite[Proposition~2.8]{EKR}
(for reduced crossed products).

\begin{prop}\label{EQOK}
Suppose $\alpha$ is an action of a discrete group $G$ on a
$C^*$-algebra $A$, and $H$ is a subgroup of $G$. Then there is
an isomorphism of the Echterhoff-Quigg crossed product
$(A\times_\alpha G)\times_{\widehat\alpha|} (G/H)$ onto
$(A\otimes c_0(G/H))\times_{\alpha\otimes\lt}G$ which
carries each representation
\begin{equation}\label{maponreps}
(\psi\times W)\times\nu\ \mbox{ into }\ (\psi\otimes\nu)\times W;
\end{equation}
in particular,
for each representation $\phi\times U$ of $A\times_\alpha G$
it carries the \qregular representation
\begin{equation}\label{maponqreg}
(((\phi\times U)\otimes\lambda)\circ\widehat\alpha)\times (1\otimes
M)\ \mbox{ into }\   (\phi\otimes M)\times (U\otimes \lambda).
\end{equation}
\end{prop}

\begin{remark}\label{otimesconv}
Equations~\eqref{maponreps} and~\eqref{maponqreg} involve tensor product
symbols with at least three different meanings.
While in most cases the meaning is clear from context, there is
potential for confusion
when trying to compare the maximal tensor product $\psi\otimes\nu$ in
\eqref{maponreps} with the spatial tensor product
$\phi\otimes M$ in \eqref{maponqreg}; this can be resolved by writing
$\phi\otimes M$ as $(\phi\otimes 1)\otimes_{\max}(1\otimes M)$.
\end{remark}

\begin{proof}[Proof of Proposition~\ref{EQOK}]
For $x\in G$,
the spectral subspace $B_x$ in the crossed  product
$B=A\times_\alpha G$ is given in terms of the universal
covariant representation $(i_A,i_G)\colon (A,G)\to
M(A\times_\alpha G)$ by $B_x:=\{i_A(a)i_G(x)\mid a\in A\}$.
Thus (using \eqref{defn-cov-rep})
a triple $(\psi,W,\nu)$ gives a covariant representation
$(\psi\times W,\nu)$ of $(A\times_\alpha G,
G/H,\widehat\alpha|)$ if and only if
\[
\psi(a)W_x\nu(\epsilon_{yH})=\nu(\epsilon_{xyH})\psi(a)W_x\ \mbox{
for all $x,y\in G$ and $a\in A$,}
\]
which happens if and only if the ranges of $\psi$ and $\nu$ commute
and $(\nu,W)$ is a  covariant representation of
$(c_0(G/H),G,\lt)$. In particular, this observation implies
that the canonical embeddings $(k_A,k_G, k_{c(G/H)})$ of
$(A,G,c_0(G/H))$ in $M((A\otimes
c_0(G/H))\times_{\alpha\otimes\lt}G)$ generate a covariant
representation $(k_A\times k_G,k_{c(G/H)})$ of
$(A\times_\alpha G, G/H,\widehat\alpha|)$, and hence by
\cite[Proposition~2.7]{EQ} induce a homomorphism $\Lambda$ of
$(A\times_\alpha G)\times_{\widehat\alpha|} (G/H)$ into
$M((A\otimes c_0(G/H))\times_{\alpha\otimes\lt}G)$ such that
\begin{equation}\label{covtriple}
\Lambda(i_A(a)i_G(x),yH)
=k_A(a)k_G(x)k_{c(G/H)}(\epsilon_{yH}).
\end{equation}
We shall prove that $\Lambda$ is the required isomorphism.

First observe that elements of the form~\eqref{covtriple}
belong to and span a dense subspace of
$(A\otimes c_0(G/H))\times_{\alpha\otimes\lt}G$,
so $\Lambda$ is surjective.
To see that $\Lambda$ is injective, we show that every
representation of $(A\times_\alpha G)\times_{\widehat\alpha}
(G/H)$ factors through $\Lambda$.
Since $\widehat\alpha$ is maximal, we know from
\cite[Proposition~2.7]{EQ} that every such representation has the
form $\theta\times \nu$ for some covariant pair $(\theta,\nu)$, and
we can write $\theta$ as $\psi\times W$ for some covariant
representation $(\psi,W)$ of $(A,G,\alpha)$. The description of the
covariant representations in the previous paragraph implies, first,
that $\psi$ and $\nu$ combine to
give a representation $\psi\otimes \nu$ of $A\otimes c_0(G/H)$
(since their ranges commute),
and second, that
$(\psi\otimes\nu,W)$ is covariant for $\alpha\otimes \lt$
(since
$(\psi,W)$ is covariant for $\alpha$ and
$(\nu,W)$ is covariant for $\lt$).
Now we just need to check using
\eqref{covtriple} that
\begin{equation}\label{phigoodonreps}
((\psi\otimes\nu)\times W)\circ \Lambda=(\psi\times W)\times \nu,
\end{equation}
and deduce that $\Lambda$ is injective. The formula
\eqref{phigoodonreps} immediately gives \eqref{maponreps}.

To see what $\Lambda$ does to the \qregular representation
induced by a representation  $\phi\times U$ of $A\times_\alpha G$,
recall that the dual
coaction $\widehat \alpha$ satisfies
$\widehat\alpha(i_A(a))=i_A(a)\otimes 1$ and
$\widehat\alpha(i_G(x))=i_G(x)\otimes u(x)$,
where $u\colon G\to C^*(G)$ is the canonical map.
So
$((\phi\times U)\otimes\lambda)\circ\widehat\alpha
=(\phi\otimes1)\times(U\otimes\lambda)$,
and thus, according to \eqref{phigoodonreps} and Remark~\ref{otimesconv},
\[
(((\phi\times U)\otimes\lambda)\circ\widehat\alpha)\times(1\otimes M)
=((\phi\otimes1)\times(U\otimes\lambda))\times(1\otimes M)
\]
is carried into
\[
((\phi\otimes 1)\otimes_{\max}(1\otimes M))\times (U\otimes\lambda)
= (\phi\otimes M)\times (U\otimes\lambda).
\]
\end{proof}

In what follows,
we denote by $X$ the
$(A\otimes c_0(G/H))\times_{\alpha\otimes\lt}G - A\times_\alpha H$
imprimitivity bimodule constructed by Green in \cite{green-twisted},
and we use $X\dashind$ to denote the associated bijection
of $\Rep(A\times_\alpha H)$ onto
$\Rep((A\otimes c_0(G/H))\times_{\alpha\otimes\lt}G)$.

\begin{thm}\label{appltoext}
Suppose that $(G,H)$ is a discrete Hecke pair,
that $\alpha$ is an action  of  $G$ on a $C^*$-algebra $A$,
and $(\psi,T)$ is a covariant representation of $(A, H,\alpha)$
on some Hilbert space $\H_\psi$.
Let $\H = X\otimes_{A\times_\alpha H}\H_\psi$,
and let $\pi\colon A\times_\alpha G\to B(\H)$
and $\nu\colon c_0(G/H)\to B(\H)$
be such that
$\pi\times\nu$ is the representation of
$(A\times_\alpha G)\times_{\widehat\alpha|}(G/H)$
corresponding to $X\dashind(\psi\times T)$
under the isomorphism of Proposition~\ref{EQOK}.
Then
there exists a representation $(\psi,\overline{T})$ of $(A, G,\alpha)$ on
$\H_\psi$ such that $\overline{T}|_H=T$
if and only if
there exists a representation $V$
of $\H(G,H)$ on $\H$ such that
$(\nu,V)$ is a covariant pair and such that the
range of $V$ commutes with the range of $\pi$.

In the special case where
$G=N\rtimes Q$, $S$, and $H$ are as in
Section~\ref{two} and Corollary~\ref{cor-ext-LL},
there exists such a representation $(\psi,\overline{T})$
if and only if there
exist a unitary representation
$U\colon N/H\to B(\H)$
and an isometric representation
$W\colon S\to B(\H)$ which, together with $\pi$ and $\nu$,
satisfy conditions~\textup{(i)--(iv)} of Corollary~\ref{cor-ext-LL}.
\end{thm}

\begin{proof}
By \cite[Theorem~1]{aHKR-limbo}, there exists such a
representation $(\psi,\overline{T})$ if and only if there exists
a
representation $\phi\times U$ of $A\times_\alpha G$ such that
$X\dashind(\psi\times T)$ is unitarily equivalent to
$(\phi\otimes M)\times(U\otimes\lambda)$;
by Proposition~\ref{EQOK} these are precisely
the representations $\phi\times U$ such that
$\pi\times\nu$ is unitarily equivalent to
the \qregular representation
$(((\phi\times U)\otimes\lambda)\circ\widehat\alpha)
\times(1\otimes M)$.
Thus the results follow by applying Theorem~\ref{disc-thm}
and Corollary~\ref{cor-ext-LL} to the maximal coaction
$\delta=\widehat\alpha$ of $G$ on $B=A\times_\alpha G$.
\end{proof}

In principle, whenever
$H$ is normal in $N$ and $N$ is normal in $G$
(as occurs in the second part of Theorem~\ref{appltoext}),
the question of extending representations from $H$
to $G$ is answered by applying
\cite[Theorem~4]{aHKR-limbo} twice:
a representation of $H$ extends to $G$ if and only if it extends to $N$
and the extension in turn extends to $G$.
But in practice,
the criteria for extending from $N$ to~$G$ cannot be verified
because the extension from $H$ to~$N$, when it exists,
is not explicitly constructed in \cite{aHKR-limbo}.


\end{document}